\documentclass[12pt]{article}

\usepackage{amsmath, epsfig, cite,colordvi}
\usepackage{amssymb}
\usepackage{amsfonts}
\usepackage{latexsym}
\usepackage{graphicx}
\newtheorem{thm}{Theorem}[section]

\newtheorem{conj}[thm]{Conjecture}

\newcommand{\qed}{{\hfill\rule{4pt}{7pt}}}
\def\pf{\noindent {\it Proof.} }

\numberwithin{equation}{section}

\makeatletter \@addtoreset{equation}{section} \makeatother

\begin{document}

\begin{center}
{\bf \Large Decomposition of Triply Rooted Trees}
\end{center}

\begin{center}
William Y. C. Chen$^{1}$, Janet F.F. Peng$^{2}$ and Harold R.L. Yang$^3$\\[6pt]

$^{1,2,3}$Center for Combinatorics, LPMC-TJKLC\\
Nankai University, Tianjin 300071, P. R. China

Email: $^{1}${\tt
chen@nankai.edu.cn}, $^{2}${\tt janet@mail.nankai.edu.cn}, $^3${\tt yangruilong@mail.nankai.edu.cn}
\end{center}

\vskip 6mm \noindent {\bf Abstract.} In this paper, we give a
decomposition of triply rooted trees into three doubly rooted trees.
This leads to a combinatorial interpretation of an identity
conjectured by Lacasse in the study of the PAC-Bayesian machine learning theory,
and proved by Younsi by using the Hurwitz identity on multivariate
Abel polynomials. We also give
a bijection between the set of functions from $[n+1]$ to $[n]$ and the set of triply rooted trees on $[n]$,
which leads to the refined enumeration of functions from $[n+1]$ to $[n]$
with respect to the number of elements in the orbit of $n+1$ and the number of periodic points.

\noindent{\bf Keywords:} doubly rooted tree, triply rooted tree, bijection

\noindent{\bf AMS Classification:} 05A15, 05A19

\section{Introduction}

Lacasse \cite{Lacasse} introduced the functions $\xi(n)$ and $\xi_2(n)$ in his study
of the classical PAC-Bayes theorem in the theory of machine learning, where
\[
    \xi(n)=\sum_{k=0}^n{n\choose k}\left({k\over n}\right)^k\left(1-{k\over n}\right)^{n-k}
\]
and
\[
    \xi_2(n)=\sum_{j=0}^n\sum_{k=0}^{n-j}{n\choose j}{ n-j \choose k}\left(j\over n\right)^j\left(k\over n\right)^k\left(1-{{j}\over n}-{k\over n}\right)^{n-j-k}.
\]
He showed that $\xi(n)$ can be used to give a tighter bound of the Kullback-Leibler divergence between the risk and the empirical risk on a sample space $S$ of a hypothesis function in a hypothesis space, whereas $\xi_2(n)$ can be used to bound the Kullback-Leibler divergence between the risk and the empirical risk on $S$ of the joint distribution of two hypothesis functions in a hypothesis space.

While $\xi_2(n)$ is a double sum,  based on numerical evidence Lacasse \cite{Lacasse}
posed the following conjecture stating that $\xi_2(n)$ can be reduced to the single sum $\xi(n)$.

\begin{conj}\label{Lacasse}
For $n\in \mathbb{N}$, we have
\begin{equation}\label{eqi1}
    \xi_2(n)=\xi(n)+n.
\end{equation}
\end{conj}

By applying an identity of Hurwitz on  multivariate
Abel polynomials, Younsi \cite{Youns} gave an algebraic proof of this conjecture.
Recall that multivariate Abel polynomials are defined by
\[
A_n(x_1,x_2,\ldots,x_m;p_1,p_2,\ldots,p_m)=\sum_{k_{1}+k_{2}+\cdots+k_m=n}{n\choose k_1,k_2,\ldots,k_m}\prod_{j=1}^m (x_j+k_j)^{k_j+p_j},
\]
where $x_1,x_2,\ldots,x_m\in \mathbb{R}$ and $p_1,p_2,\ldots,p_m\in \mathbb{Z}$.
Hurwitz proved that under certain conditions, the polynomials $A_n(x_1,x_2,\ldots,x_m;p_1,p_2,\ldots,p_m)$
reduce to single sums. In particular, when $p_1=p_2=\cdots=p_m=0$, we have
\begin{equation}\label{eqan}
A_n(x_1,\ldots,x_m;0,\ldots,0)=\sum_{k=0}^n{n\choose k}(x_1+x_2+\cdots+x_m+n)^{n-k}\alpha_k(m-1),
\end{equation}
where $\alpha_k(r)=r(r+1)\cdots (r+k-1)$ is the rising factorial, see, for example, Riordan \cite{Riordan}.

Younsi \cite{Youns} observed that $\xi(n)=A_n(0,0;0,0)$ and $\xi_2(n)=A_n(0,0,0;0,0,0)$, and obtained the following expressions for $\xi(n)$ and $\xi_2(n)$ by the above identity (\ref{eqan}),
\begin{eqnarray}
\xi(n)&=&{1\over n^n}\sum_{j=0}^n n^j {n!\over j!},\label{xi}\\
\xi_2(n)&=&{1\over n^n}\sum_{j=0}^n n^{n-j}{n \choose j}(j+1)!. \label{xi2}
\end{eqnarray}
 Conjecture \ref{Lacasse} can be easily deduced from (\ref{xi}) and (\ref{xi2}).

In this paper, we give a combinatorial explanation of  relation (\ref{eqi1}).
Rewriting (\ref{eqi1}) as
\begin{equation}
    \sum_{j=0}^n\sum_{k=0}^{n-j}{n\choose j}{ n-j \choose k}j^j k^k(n-j-k)^{n-j-k}=\sum_{k=0}^n{n\choose k}k^k (n-k)^{n-k}+n^{n+1},
\end{equation}
we see that it is equivalent to the following form
\begin{equation}\label{eq1}
    \sum_{j=1}^n\sum_{k=0}^{n-j}{n\choose j}{ n-j \choose k}j^j k^k(n-j-k)^{n-j-k}=n^{n+1}.
\end{equation}
The right hand side of (\ref{eq1}) indicates that
we need the notion of triply rooted trees, namely,  labeled trees with three distinguished,
but not necessarily distinct  vertices. To be more specific,
the three distinguished vertices of a triply rooted tree
are called the first, the second and the third root, respectively.
It can be easily seen  that the summand on the left hand side of (\ref{eq1})
can be interpreted as the number of
triples of doubly rooted trees with a given number of vertices in
each doubly rooted tree. Hence relation (\ref{eq1}) can be deduced
from a decomposition of a triply rooted tree into three doubly rooted trees.

The second result of this paper is a correspondence between the set of
 functions from $[n+1]$ to $[n]$ and the set of triply rooted trees on $[n]$. Let $f$ be a function from $[n+1]$ to $[n]$ and let $T$ be the corresponding triply rooted tree. We find that the orbit of $n+1$ on $f$  is mapped to the set of ancestors of the second root in $T$, and the set of periodic points of $f$ is mapped to the set of ancestors of the third root in $T$.
Based on this property of our bijection, we derive a formula for the number of functions from $[n+1]$ to $[n]$ with a given number of elements in
  the orbit of $n+1$ and a given number of periodic points.

\begin{section}{Decomposition of triply rooted tree}\label{decomposition}

In this section, we give a combinatorial interpretation of Lacasse's identity
 by providing a decomposition of a triply rooted tree into
 three doubly rooted trees.

Recall that a rooted tree is defined to be a labeled tree with a specific vertex, which is called the root.
Let $R_n$ denote the set of rooted trees on $[n]$. The set $R_n$
 is counted by
$n^{n-1}$, see Stanley \cite[5.3.2 Proposition]{Stanley}.
A doubly rooted tree is defined as a labeled tree with two distinguished vertices $r_1$ and $r_2$, where we call $r_1$ the first root and call $r_2$ the second root.
 Notice that the two roots of a doubly rooted tree are not
 required to be distinct. We denote by $D_n$ the set of doubly rooted trees on $[n]$. From the formula for $|R_n|$, one sees that $|D_n|=n^n$. The notion of
 doubly rooted trees leads to be an elegant proof of the formula for $|D_n|$
 independently obtained  by Goulden and Jackson \cite{GJ} and
Joyal \cite{Joyal}.

The identity (\ref{eq1}) indicates that there is a decomposition of a
triply rooted tree into three doubly rooted trees. More precisely, we define a triply rooted tree to be a labeled tree with three distinguished vertices $ r_1,r_2$, and $r_3 $, which are called the first, the second, and the third root, respectively. Again, the three roots of a triply rooted tree are not necessarily distinct. Denote by $T_n$ the set of triply rooted trees.  From the formula for $|D_n|$ it is clear that $|T_n|=n^{n+1}$. So the right hand side of (\ref{eq1})
can be interpreted as the number of triply rooted trees on $[n]$.

On the other hand, let $Q_n$ denote the set of triples of doubly rooted trees $(D,D',D'')$ such that the vertex sets of $D$, $D'$, $D''$ form a composition of $[n]$ with $D$ being nonempty. To be more specific, a triple
$(X,Y,Z)$ of subsets of a set $S$ is said to be a composition
 of $S$ if $X,Y$, and $Z$ are disjoint and their union equals $S$.
 It is obvious that $Q_n$ is counted by
\[
\sum_{j=1}^n\sum_{k=0}^{n-j}{n\choose j}{ n-j \choose k}j^j k^k(n-j-k)^{n-j-k},
\]
which is the left hand side of (\ref{eq1}).
Hence identity (\ref{eq1}) follows from the following bijection.

\begin{thm}\label{thm1}
For $n\geq 1$, there is a bijection between $Q_n$ and $T_n$.
\end{thm}

To present the proof of the above theorem, we recall some terminology.
Given two vertices $i$ and $j$ of a rooted tree $T$, we say that $j$
is a \emph{descendant} of $i$, or $i$ is an \emph{ancestor} of $j$,
 if $i$ lies on the unique path from
the root to $j$.
In particular, each vertex is a descendant as well as an ancestor of
itself. A \emph{child} of $i$ means a descendant $j$ of $i$ such that $(i,j)$ is an edge of $T$.
The depth of $i$ is defined to be the number of edges of the unique path from the root to $i$.
Given two vertices $v_1$ and $v_2$ of $T$, there is a unique
vertex $v$  that is the common ancestor of $v_1$ and $v_2$ with the largest depth.
This vertex is called the \emph{least common ancestor} of $v_1$ and $v_2$, see Aho, Hopcroft and Ullman \cite{AHU}.
For example, for the tree in Figure \ref{fig1}, the least common ancestor of  $1$ and   $3$ is  $5$, while the least common ancestor of   $1$ and $6$ is the root $4$.

\begin{figure}[h,t]
\begin{center}
\begin{picture}(100,60)
\setlength{\unitlength}{1mm}

\put(16,15){\circle*{2}}\put(15.5,17){\small $4$}
\put(10,5){\circle*{2}}\put(7,4.5){\small $2$}
\put(10,5){\line(3,5){6}}

\put(22,5){\circle*{2}}\put(23.5,4.5){\small $5$}
\put(22,5){\line(-3,5){6}}

\put(10,-5){\circle*{2}}\put(9.5,-9){\small $6$}
\put(10,-5){\line(0,1){10}}

\put(16,-5){\circle*{2}}\put(15.5,-9){\small $3$}
\put(16,-5){\line(3,5){6}}

\put(28,-5){\circle*{2}}\put(27.5,-9){\small $1$}
\put(28,-5){\line(-3,5){6}}
\end{picture}
\vspace{1cm}\caption{A rooted tree on $[6]$.}\label{fig1}
\end{center}
\end{figure}
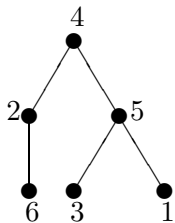

Throughout this paper, we use
$r_1(D)$ and $r_2(D)$ to denote the first root and the second root of
a doubly rooted tree $D$, respectively, and we use $r_1(T)$, $r_2(T)$ and $r_3(T)$ to denote
the first root, the second root and the third root of a triply rooted tree $T$, respectively.

\noindent
{\it Proof of Theorem \ref{thm1}.}
We define a map $\varphi$ from $Q_n$ to $T_n$. Given a triple $(D,D',D'')$
 of doubly rooted trees in $Q_n$, we aim to construct
 a triply rooted tree on $[n]$. First, we consider the case when neither $D'$ nor $D''$ is empty.

 We merge $D$ and $D'$ by
 setting $r_1(D')$ to be a child of $r_2(D)$,
 and we merge $D$ and  $D''$ by setting $r_1(D'')$ to be a child of $r_2(D)$.
 By setting $r_2(D')$ and $r_2(D'')$ to be the second root and the third root of the resulting tree, we obtain a triply rooted tree $T$.

For example, Figure \ref{fig2} gives an illustration of a triple of doubly rooted trees and the corresponding triply rooted tree, where the second
root is represented by a solid square, and the third root is represented by a hollow square.

\begin{figure}[h,t]
\begin{center}
\begin{picture}(410,80)
\setlength{\unitlength}{1mm}

\put(16,15){\circle*{2}}\put(15.5,17){\small $6$}
\put(9,4){{\hfill \rule{6pt}{6pt}}}\put(5,4.5){\small $12$}
\put(10,5){\line(3,5){6}}

\put(22,5){\circle*{2}}\put(23.5,4.5){\small $8$}
\put(22,5){\line(-3,5){6}}

\put(10,-5){\circle*{2}}\put(9.5,-9){\small $9$}
\put(10,-5){\line(0,1){10}}

\put(16,-5){\circle*{2}}\put(15.5,-9){\small $3$}
\put(16,-5){\line(3,5){6}}

\put(28,-5){\circle*{2}}\put(27.5,-9){\small $1$}
\put(28,-5){\line(-3,5){6}}

\put(35,14){\hfill \rule{6pt}{6pt}}\put(35,17){\small $2$}
\put(36,5){\line(0,1){10}}\put(36,5){\circle*{2}}\put(34,1){\small $10$}

\put(54,15){\circle*{2}}\put(53.5,17){\small $5$}
\put(54,5){\line(0,1){10}}\put(54,5){\circle*{2}}
\put(53,4){\hfill \rule{6pt}{6pt}}\put(51,4.5){\small $4$}
\put(48,-5){\circle*{2}}\put(48,-5){\line(3,5){6}}\put(47,-9){\small $7$}
\put(60,-5){\circle*{2}}\put(60,-5){\line(-3,5){6}}\put(58,-9){\small $11$}
\put(14,-15){$D$}\put(34,-15){$D'$}\put(52,-15){$D''$}

\put(80,5){$\Longrightarrow$}

\put(0,10){\begin{picture}(0,0)
\put(120,15){\circle*{2}}\put(119.5,17){\small $6$}
\put(110,5){\circle*{2}}\put(110,5){\line(1,1){10}}
\put(105,4.5){\small $12$}
\put(110,-5){\circle*{2}}\put(110,-5){\line(0,1){10}}
\put(109.5,-9){\small $9$}
\put(116,-5){\circle*{2}}\put(116,-5){\line(-3,5){6}}
\put(117.5,-5.5){\small $5$}
\put(115,-16){\line(1,0){2}}\put(115,-14){\line(1,0){2}}
\put(115,-16){\line(0,1){2}}\put(117,-16){\line(0,1){2}}
\put(116,-14){\line(0,1){8}}
\put(117.5,-15.5){\small $4$}
\put(110,-25){\circle*{2}}\put(110,-25){\line(3,5){5.35}}
\put(109,-29){\small $7$}
\put(122,-25){\circle*{2}}\put(122,-25){\line(-3,5){5.35}}
\put(120.5,-29){\small $11$}

\put(103,-6){\hfill \rule{6pt}{6pt}}
\put(104,-5){\line(3,5){6}}
\put(100.5,-5.5){\small $2$}\put(104,-15){\circle*{2}}
\put(104,-15){\line(0,1){10}}\put(102,-19){\small $10$}

\put(130,5){\circle*{2}}\put(130,5){\line(-1,1){10}}
\put(131.5,4.5){\small $8$}
\put(124,-5){\circle*{2}}\put(124,-5){\line(3,5){6}}
\put(123,-9){\small $3$}
\put(136,-5){\circle*{2}}\put(136,-5){\line(-3,5){6}}
\put(135,-9){\small $1$}

\end{picture}
}
\end{picture}
\vspace{2cm}\caption{The merging process when $D'\not=\emptyset$ and $D''\not=\emptyset$.}\label{fig2}
\end{center}
\end{figure}
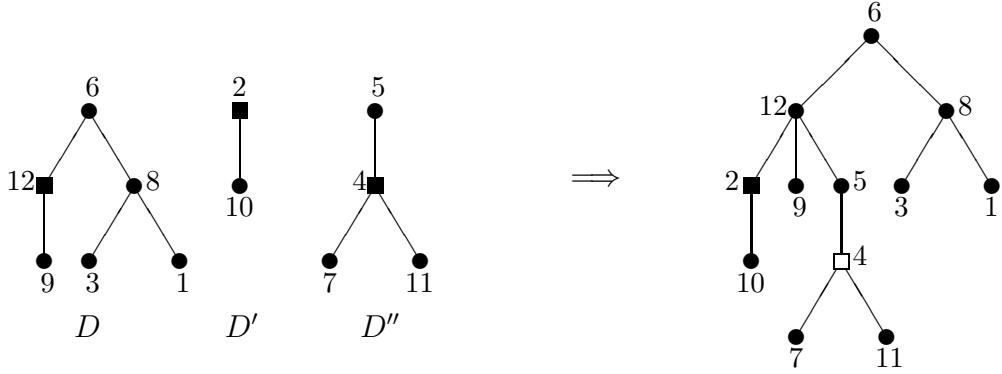

We now consider the case when either $D'$ or $D''$ is empty. If $D'=\emptyset$ and $D''\neq \emptyset$, we merge $D$ and $D''$ by setting  $r_1(D'')$ to be a child of $r_2(D)$. Setting $r_2(D)$ and $r_2(D'')$ to be the second root and the third root, we  obtain a triply rooted tree $T$.

 If $D'\neq \emptyset$ and $D''=\emptyset$, we merge $D$ and $D'$ by setting  $r_1(D')$ to be a child of $r_2(D)$. Setting $r_2(D')$ and $r_2(D)$ to be the second root and the third root, we obtain a triply rooted tree $T$.

 If both $D$ and $D''$ are empty, then we set $r_2(D)$ to be the second root and the third root to obtain a triply rooted tree $T$.

In summary,  $(r_1(T),r_2(T),r_3(T))$  is given as follows:
\[\left\{
\begin{array}{ll}
(r_1(D),r_2(D'),r_2(D'')),&  \mbox{\ \ if $D'\neq\emptyset$, and $D''\neq\emptyset$};\\[3pt]
(r_1(D),r_2(D),r_2(D'')),&\mbox {\ \ if $D'=\emptyset$, and $D''\neq\emptyset$};\\[3pt]
(r_1(D),r_2(D'),r_2(D)),& \mbox {\ \ if $D'\neq\emptyset$, and $D''=\emptyset$};\\[3pt]
(r_1(D),r_2(D),r_2(D)),& \mbox {\ \ if $D'=\emptyset$, and $D''=\emptyset$}.
\end{array}\right.
\]
Figure \ref{fig3} gives an illustration of the
 merging process when $D'=\emptyset$ and $D''\neq \emptyset$.

\begin{figure}[h,t]
\begin{center}
\begin{picture}(410,80)
\setlength{\unitlength}{1mm}

\put(16,15){\circle*{2}}\put(15.5,17){\small $6$}
\put(9,4){\hfill \rule{6pt}{6pt}}\put(5,4.5){\small $10$}
\put(10,5){\line(3,5){6}}

\put(22,5){\circle*{2}}\put(23.5,4.5){\small $8$}
\put(22,5){\line(-3,5){6}}

\put(10,-5){\circle*{2}}\put(9.5,-9){\small $9$}
\put(10,-5){\line(0,1){10}}

\put(16,-5){\circle*{2}}\put(15.5,-9){\small $3$}
\put(16,-5){\line(3,5){6}}

\put(28,-5){\circle*{2}}\put(27.5,-9){\small $1$}
\put(28,-5){\line(-3,5){6}}

\put(54,15){\circle*{2}}\put(53.5,17){\small $5$}
\put(54,5){\line(0,1){10}}\put(53,4){\hfill \rule{6pt}{6pt}}
\put(51,4.5){\small $4$}
\put(48,-5){\circle*{2}}\put(48,-5){\line(3,5){6}}\put(47,-9){\small $7$}
\put(60,-5){\circle*{2}}\put(60,-5){\line(-3,5){6}}\put(59,-9){\small $2$}
\put(14,-15){$D$}\put(34,-15){$D'$}\put(35,5){$\emptyset$}
\put(52,-15){$D''$}

\put(80,5){$\Longrightarrow$}

\put(0,10){\begin{picture}(0,0)
\put(120,15){\circle*{2}}\put(119.5,17){\small $6$}
\put(109,4){\hfill \rule{6pt}{6pt}}\put(110,5){\line(1,1){10}}
\put(105,4.5){\small $10$}

\put(116,-5){\circle*{2}}\put(116,-5){\line(-3,5){6}}
\put(117.5,-5.5){\small $5$}
\put(115,-16){\line(1,0){2}}\put(115,-14){\line(1,0){2}}
\put(115,-16){\line(0,1){2}}\put(117,-16){\line(0,1){2}}
\put(116,-14){\line(0,1){8}}
\put(117.5,-15.5){\small $4$}
\put(110,-25){\circle*{2}}\put(110,-25){\line(3,5){5.4}}
\put(109,-29){\small $7$}
\put(122,-25){\circle*{2}}\put(122,-25){\line(-3,5){5.4}}
\put(121,-29){\small $2$}

\put(104,-5){\circle*{2}}
\put(104,-5){\line(3,5){6}}
\put(100.5,-5.5){\small $9$}

\put(130,5){\circle*{2}}\put(130,5){\line(-1,1){10}}
\put(131.5,4.5){\small $8$}
\put(124,-5){\circle*{2}}\put(124,-5){\line(3,5){6}}
\put(123,-9){\small $3$}
\put(136,-5){\circle*{2}}\put(136,-5){\line(-3,5){6}}
\put(135,-9){\small $1$}

\end{picture}
}
\end{picture}
\vspace{2cm}\caption{The merging process when $D'=\emptyset$ and $D''\not=\emptyset$.}\label{fig3}
\end{center}
\end{figure}
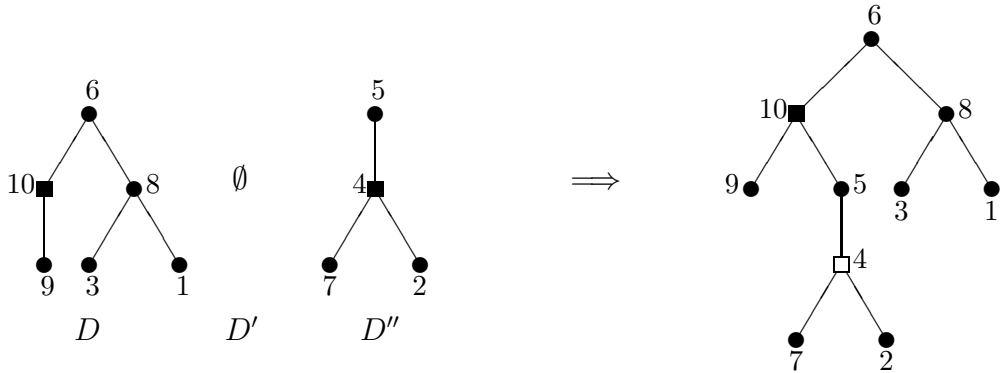

To show that the above process is invertible, we give a description of the inverse
procedure. Given a
triply rooted tree $T$ with three roots $r_1,r_2$ and $r_3$, assume that
$w$ is the least
common ancestor of $r_2$ and $r_3$.

 We first consider the
case when $w\neq r_2$ and $w\neq r_3$. We proceed to find two edges such that
by removing them we can recover three doubly rooted trees
$D$, $D'$ and $D''$.
Find the child $x$  of $w$ such that $r_2$ is  a descendant of $x$,
and the child $y$ of $w$ such that $r_3$ is a descendant of $y$.
Removing the edges $(w,x)$ and $(w,y)$, we  get three
trees with three roots $r_1$, $x$ and $y$. Let $D$ be the doubly rooted trees with two roots $r_1$ and $w$, let $D'$ be the doubly rooted trees with two roots $x$ and $r_2$, and
let $D''$ be the doubly rooted trees with two roots $y$ and $r_3$.

Next we consider the remaining cases.
When $w=r_2$ and $w\not=r_3$, suppose that $y$ is the child of $w$ such that $r_3$ is a descendant of $y$. Removing the edge $(w,y)$, we get
 a doubly rooted tree $D$ with two roots $r_1$ and $w$,
 a doubly rooted tree
 $D''$  with two roots $y$ and $r_3$. Moreover, we set $D'=\emptyset$.

When $w=r_3$ and $w\not=r_2$, suppose that $x$ is the child of $w$ such that $r_2$ is a descendant of $x$. Removing the edge $(w,x)$, we get a doubly rooted tree $D$ with two roots $r_1$ and $w$ and a doubly rooted tree $D'$ with two roots $x$ and $r_2$.
 Meanwhile, we set $D''=\emptyset$.

When $w=r_2=r_3$, let $D$ to be the doubly rooted tree obtained
from $T$ by setting the two roots to be $r_1$ and $w$, and let
 $D'=\emptyset$ and $D''=\emptyset$.

It can be easily checked that in any case
the three doubly rooted trees $D$, $D'$ and $D''$ can be merged
into the triply rooted tree $T$. That is, the above merging process
is invertible.
This completes the proof. \qed

\end{section}

\section{Functions from $[n+1]$ to $[n]$}\label{triply}

In this section, we establish a correspondence between
functions from
$[n+1]$ to $[n]$ and triply rooted trees on $[n]$,
 which maps the orbit of $n+1$ to
 the set of ancestors of the second root, and maps the set of periodic points to the set of ancestors of the
 third root. By the symmetry between the
 second and third roots, we deduce a symmetry property
 of the number of functions from
 $[n+1]$ to $[n]$ with respect to the number of periodic points and
 the size of the orbit of $n+1$.

Given a function $f$ from $[n+1]$ to $[n]$, the \emph{orbit} of $x$ on $f$ is defined to be the set $\{x, f(x),f^2(x),\ldots\}$. If there exists some $j\geq 1$, such that $f^j(x)=x$, then $x$ is called a \emph{periodic point}
of $f$. We have the following
correspondence.

\begin{thm}\label{thm2}
  There is a bijection $\phi$ between the set of functions  $f$ from $[n+1]$ to $[n]$ and the set of triply rooted trees  on $[n]$ such that
     the orbit of $n+1$ on $f$ excluding $n+1$ itself is mapped
       to the set of ancestors of the second root of $\phi(f)$ and the set of periodic points of $f$ is mapped to the set of ancestors of the third root
       of $\phi(f)$.
\end{thm}

\pf The map $\phi$  can be described as follows.
  Let $f$ be a function from $[n+1]$ to $[n]$.
  We proceed to construct a triply rooted tree $T$ on $[n]$ based on the function
  $f$.  We begin with the functional digraph $G_f$  of $f$, that is,
  a digraph on $[n+1]$ with arcs $(i, f(i))$ for $1\leq i \leq n+1$.
  Let $C_1$ be the connected component of $G_f$ containing the vertex $n+1$.
  Consider the longest path $P$ starting from $n+1$, say,
\[P\colon n+1=u_0\rightarrow u_1\rightarrow u_2
\rightarrow \cdots \rightarrow u_k.\]
In other words, $k$ is the smallest integer such that
$f(u_k)=u_j$ for some $j\leq k$. Removing the arc $(u_k,u_{j})$ and the vertex $n+1$ from $C_1$, we get a tree $H$ rooted at $u_k$.

Let $C_2$ be the digraph $G_f\setminus C_1$. When $C_2=\emptyset$,  we set $u_k,u_1$ and $u_j$ to be the three roots of $H$ to obtain a triply rooted tree $T$.

When $C_2\not=\emptyset$, suppose that the vertex set of $C_2$ is $\{v_1,v_2,\ldots,v_s\}$. Note that $C_2$ is a functional digraph on
 $\{v_1, v_2, \ldots, v_s\}$. By applying the bijection
 between functions and doubly rooted trees, obtained by
 Joyal \cite{Joyal} and Goulden and Jackson \cite{GJ}, $C_2$ corresponds to a doubly rooted tree $D$ on $\{v_1,v_2,\ldots,v_s\}$. Let $w_1$ and $w_2$ be the two roots of $D$.

Finally, we merge the rooted tree $H$ and
 the doubly rooted tree $D$ by joining the first root $w_1$ of $D$
  and the vertex $u_j$ of $H$  with $w_1$ being the child.
    Setting $u_k$, $u_1$ and $w_2$ to be the first, the second and the third root, respectively, we get a triply rooted tree $T$, and we set $\phi(f)=T$.

\begin{figure}[h,t]
\begin{center}
\begin{picture}(300,50)
\setlength{\unitlength}{1mm}

\multiput(10,10)(10,0){5}{\circle*{1.5}}
\put(10,10){\line(1,0){40}}
\multiput(10,10)(10,0){3}{\vector(1,0){6}}\put(50,10){\vector(-1,0){6}}
\put(20,0){\circle*{1.5}}
\put(20,0){\vector(0,1){6}}\put(20,0){\line(0,1){10}}
\put(40,0){\circle*{1.5}}\put(40,0){\line(0,1){10}}
\put(40,10){\vector(0,-1){6}}
\put(40,0){\line(1,0){10}}
\put(40,0){\vector(-1,1){6}}
\put(40,0){\line(-1,1){10}}
\put(50,0){\circle*{1.5}}\put(50,0){\vector(-1,0){6}}
\multiput(80,10)(10,0){2}{\circle*{1.5}}
\multiput(70,0)(10,0){2}{\circle*{1.5}}\put(70,0){\line(1,1){10}}
\put(80,0){\line(0,1){10}}
\put(70,0){\vector(1,1){6}}\put(80,0){\vector(0,1){6}}
\qbezier(80,10)(85,6)(90,10)\qbezier(80,10)(85,14)(90,10)
\put(85.2,11.8){\vector(1,0){1}}\put(85.5,8.1){\vector(-1,0){1}}
\put(9.2,12){\small $1$}\put(19.2,12){\small $8$}\put(29.2,12){\small $6$}\put(38.5,12){\small $12$}\put(49.2,12){\small $9$}
\put(19.2,-4){\small $3$}\put(39.2,-4){\small $2$}\put(48.5,-4){\small $10$}
\put(77,9){\small $4$}\put(91,9){\small $5$}
\put(69.3,-4){\small $7$}\put(78.5,-4){\small $11$}
\put(10,0){\vector(1,0){6}}\put(16,0){\line(1,0){4}}
\put(10,0){\circle*{1.5}}
\put(8,-4){\small $13$}
\end{picture}
\vspace{0.5cm}\caption{The functinal digraph $G_f$.}\label{fig5}
\end{center}
\end{figure}
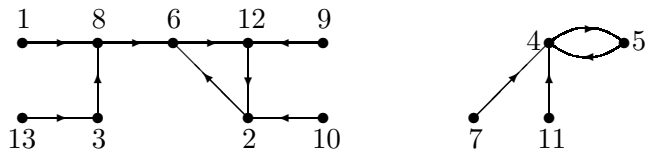

\begin{figure}[h,t]
\begin{center}
\begin{picture}(405,80)
\setlength{\unitlength}{1mm}

\put(20,15){\circle*{2}}\put(19.5,17){\small $2$}
\put(14,5){\circle*{2}}\put(14,5){\line(3,5){6}}
\put(8.5,4){\small $10$}
\put(26,5){\circle*{2}}\put(26,5){\line(-3,5){6}}
\put(27.5,4.5){\small $12$}
\put(21.2,-3){\line(3,5){5}}
\put(22.5,-5){\small $6$}
\put(20,-5){\line(1,0){2}}\put(20,-3){\line(1,0){2}}
\put(20,-5){\line(0,1){2}}\put(22,-5){\line(0,1){2}}
\put(31.4,-4){\circle*{2}}\put(31.4,-4){\line(-3,5){5}}
\put(32,-3){\small $9$}
\put(21,-12){\circle*{2}}\put(21,-12){\line(0,1){6.8}}
\put(22,-13){\small $8$}
\put(16.2,-20){\circle*{2}}\put(16.2,-20){\line(3,5){5}}
\put(13,-20){\small $1$}
\put(24.8,-21){\hfill \rule{6pt}{6pt}}\put(25.8,-20){\line(-3,5){5}}
\put(27,-20){\small $3$}

\put(14,-35){Tree $H$}

\put(60,15){\circle*{2}}\put(59.5,17){\small $5$}
\put(59,4){\hfill \rule{6pt}{6pt}}\put(60,5){\line(0,1){10}}
\put(61.5,4.5){\small $4$}
\put(54,-5){\circle*{2}}\put(54,-5){\line(3,5){6}}
\put(50.5,-5.5){\small $7$}
\put(66,-5){\circle*{2}}\put(66,-5){\line(-3,5){6}}
\put(68,-5.5){\small $11$}
\put(54,-35){Tree $D$}

\put(80,5){$\Longrightarrow$}

\put(90,10)
{
\begin{picture}(40,40)
\put(20,15){\circle*{2}}\put(19.5,17){\small $2$}
\put(10,5){\circle*{2}}\put(10,5){\line(1,1){10}}
\put(5,4.5){\small $10$}
\put(30,5){\circle*{2}}\put(30,5){\line(-1,1){10}}
\put(31.5,4.5){\small $12$}
\put(21,-4){\circle*{2}}\put(21,-4){\line(1,1){9}}
\put(23,-5){\small $6$}
\put(39,-4){\circle*{2}}\put(39,-4){\line(-1,1){9}}
\put(36,-7){\small $9$}
\put(13,-16){\circle*{2}}\put(13,-16){\line(2,3){8}}
\put(10,-16){\small $8$}
\put(7,-26){\circle*{2}}\put(7,-26){\line(3,5){6}}
\put(6,-30){\small $1$}
\put(19,-26){\circle*{2}}\put(19,-26){\line(-3,5){6}}
\put(18,-30){\small $3$}

\put(29,-16){\circle*{2}}\put(29,-16){\line(-2,3){8}}
\put(31,-16){\small $5$}
\put(29,-25.2){\line(0,1){8.8}}
\put(30,-25){\small $4$}
\put(23,-35){\circle*{2}}\put(23,-35){\line(2,3){5.3}}
\put(22,-39){\small $7$}
\put(35,-35){\circle*{2}}\put(35,-35){\line(-2,3){5.3}}
\put(34,-39){\small $11$}
\put(17.8,-27){\hfill \rule{6pt}{6pt}}
\put(28,-25){\line(1,0){2}}\put(28,-27){\line(1,0){2}}
\put(28,-27){\line(0,1){2}}\put(30,-27){\line(0,1){2}}
\end{picture}
}
\put(107,-35){Tree $T$}
\end{picture}
\vspace{4cm}\caption{An example of the bijection $\phi$.}\label{fig6}
\end{center}
\end{figure}
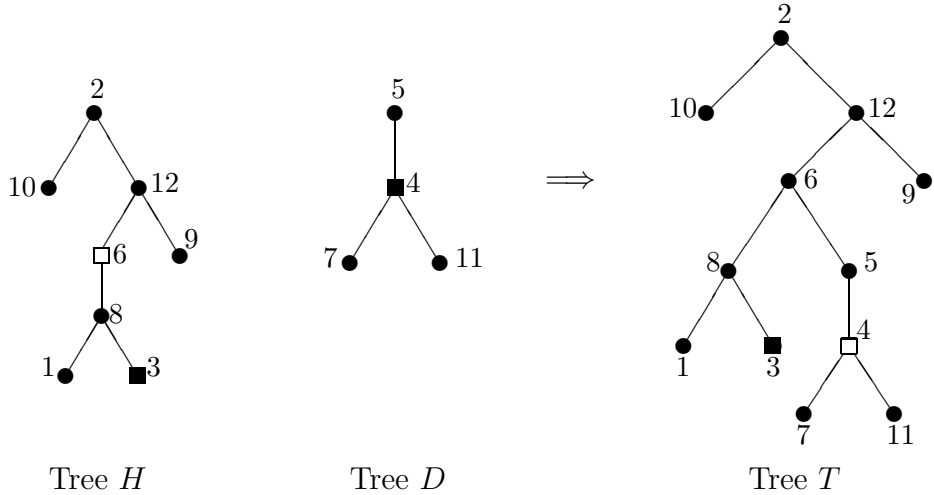

For example, let $f$ be the following function from $[13]$ to $[12]$,
\[
f=\left(
\begin{array}{ccccccccccccc}
1&2&3&4&5&6&7&8&9&10&11&12&13\\
8&6&8&5&4&12&4&6&12&2&4&2&3
\end{array}
\right).
\]
The functional digraph of $f$ is given in Figure \ref{fig5}, where $C_1$ is the
functional digraph on  $\{1,2,3,6,8,9,10,12,13\}$ and $C_2$ is the functional digraph  on $\{4,5,7,11\}$.
The longest path starting from $13$ is
\[P\colon 13\rightarrow 3\rightarrow 8\rightarrow 6\rightarrow 12 \rightarrow2,\]
with $f(2)=6$, that is, $u_1=3$, $u_k=2$ and $u_j=6$ as in the proof.
Deleting the arc $(2,6)$ and vertex $13$, we get a rooted tree $H$ as
 illustrated in Figure \ref{fig6}. By applying
  the bijection between functional digraphs and doubly rooted trees,
  $C_2$ can be mapped to a doubly rooted tree $D$ with roots $5$ and $4$
   as shown in Figure \ref{fig6}, where $w_1=5$ and $w_2=4$ as in the proof. Merging $H$ and $D$ by adding an edge $(6,5)$ and setting $2,3$ and $4$ to be the three roots, we get a triply rooted tree $T$ in Figure $\ref{fig6}$.

The map $\phi$ is indeed a bijection. The inverse map can be described as follows. For a triply rooted tree $T\in T_n$ with the three roots $r_1,r_2$ and $r_3$, we first find
the least common ancestor of
$r_2$ and $r_3$, and we denote it by $u_0$. Suppose that the unique path $P$ from $u_0$ to $r_3$ in $T$ is $u_0u_1\cdots u_i=r_3$. Removing the edge $(u_0,u_1)$ from $T$,  we get two components $T_1$ and $T_2$, where $T_1$ is rooted at $r_1$ and $T_2$ is rooted at $u_1$. Adding $n+1$ to $T_1$ by setting it as a child of $r_2$. Now, $T_1$ can be viewed as a directed graph
by making each edge point to the father. Then
 we add the arc $(r_1,u_0)$ to $T_1$ to obtain a connected
 functional digraph $C_1$.

Next, we transform $T_2$ rooted at $u_1$ into a doubly rooted tree
 by setting $r_3$ to be the second root. Then we get a
  functional digraph $C_2$ by applying the inverse map of the bijection of Joyal \cite{Joyal} and Goulden and Jackson \cite{GJ}.

Finally, let  $G=C_1\cup C_2$.  It is easily seen that $G$ is a directed graph on $[n+1]$ such that each vertex has outdegree one  and the vertex $n+1$
 has indegree zero. In other words,  $G$
  is the functional digraph of a function from $[n+1]$ to $[n]$.
 It can be checked that the above procedure is indeed the inverse of the map $\phi$.

It remains to prove the properties of $\phi$ as stated in the theorem.
 For a function $f$ from $[n+1]$ to $[n]$,   an element $x$ is a periodic point in $f$ if and only if it is a vertex in a cycle in the functional digraph $G_f$.
 It can be seen that $x$ is in a cycle
  if and only if it is an ancestor of the third root in the
   triply rooted tree $\phi(f)$. Moreover, it can be checked
    that each element $y$ in the orbit of $n+1$
   on $f$ other than $n+1$ itself corresponds to  an ancestor of the second root in the triply rooted tree $\phi(f)$.
This completes the proof.\qed

For example, for the function $f$ in Figure $\ref{fig5}$, there are five periodic points $2,12,6,5,4$, which are the vertices in the path from the root $2$ to the third root $4$ in $\phi(f)$ as demonstrated in Figure $\ref{fig6}$. The orbit of $13$ consists of $13,3,8,6,12,2$. These elements
 $3, 8,6,12,2$ correspond to the vertices in the path from the root $2$ to the second root $3$ in $\phi(f)$.

From the above bijection $\phi$, we obtain a formula for
 the number of functions from $[n+1]$ to $[n]$ with a given
   number of elements in the orbit of $n+1$ and a given number of periodic points.
This formula implies a symmetry property, which can also be interpreted
 in terms of triply rooted trees.

\begin{thm}\label{thm3} For $n\geq 1$,
let $W_{n,i,j}$ denote the set of triply rooted trees on $[n]$ such that
 the depth of the second root is $i$ and the depth
of the third root is $j$. Then we have
\begin{equation}\label{eqt}
|W_{n,i,j}|=
\sum_{d=0}^{\min(i,j)}\frac{(i+j-d+1) n!}{(n-i-j+d-1)!}n^{n-i-j+d-2},
\end{equation}
where $|W_{n,i,j}|$ is the cardinality of $W_{n,i,j}$.
\end{thm}

\pf Let $W_{n,i,j}(d)$ denote the set of triply rooted trees $T$
 in $W_{n,i,j}$ such that $d$ is the depth of the least common ancestor of the second root and the third root of $T$. We proceed to show that $W_{n,i,j}(d)$ is enumerated by the summand on the right hand side of (\ref{eqt}).

Let $T$ be a triply rooted tree  in $W_{n,i,j}(d)$.
We denote by $P_1$ the path from the first root to the second root and denote by $P_2$ the path from the first root to the third root.
It can be seen that there are exactly $k=i+j-d+1$ vertices on $P_1$ and $P_2$.
Hence the number of ways to form $P_1$ and $P_2$ equals
${n!\over (n-k)!}$. Moreover, it is known that
 there are $kn^{n-k-1}$  forests consisting of $k$ rooted trees on $[n]$
  with $k$ given roots. It follows that
   $W_{n,i,j}(d)$ is enumerated by the summand on the right hand side of (\ref{eqf}). This completes the proof.\qed

Combining Theorem \ref{thm2} and Theorem \ref{thm3}, we arrive at the following formula for the refined enumeration of functions from $[n+1]$ to $[n]$.

\begin{thm} For $n\geq 1$,
let $F_{n,i,j}$ denote the set of functions from $[n+1]$ to $[n]$ such that the size of the orbit of $n+1$ is $i$ and the number of periodic points is $j$. Then we have
\begin{equation}\label{eqf}
|F_{n,i+1,j}|=\sum_{s=0}^{\min(i,j)-1}\frac{(i+j-s-1) n!}{(n-i-j+s+1)!}n^{n-i-j+s}.
\end{equation}
\end{thm}

By the symmetry of the second roots and the third roots for $T_n$, we can conclude a symmetry relation of functions from $[n+1]$ to $[n]$ concerning the size of orbit of $n+1$ and the number of periodic points, that is,
 \begin{equation}
|F_{n,i+1,j}|=|F_{n,j+1,i}|.
\end{equation}
Notice that the above symmetry is implied by
  (\ref{eqf}).


\begin{thebibliography}{99}
\bibitem{AHU}A.V. Aho, J.E. Hopcroft and J.D. Ullman, On finding lowest common ancestors in trees, SIAM J. Computing, 5 (1), 115--132, 1976.

\bibitem{GJ}I.P. Goulden and D.M. Jackson, Combinatorial Enumeration, John Wiley, New York, 1983.

\bibitem{Joyal} A. Joyal, Une th\'{e}orie combinatoire des s\'{e}ries formelles, Adv. Math. 42 (1981), 1--82.


\bibitem{Lacasse} A. Lacasse, Bornes PAC-Bayes et algorithmes {d}'apprentissage, Ph.D. Thesis, Universite Laval, Quebec, 2010.

\bibitem{Riordan} J. Riordan, Combinatorial Identities, Robert E. Krieger Publishing Co., New York, 1968.

\bibitem{Stanley} R.P. Stanley, Enumerative Combinatorics, Vol. 2,
Cambridge University Press, Cambridge, 1999.

\bibitem{Youns} M. Younsi, Proof of a combinatorial conjecture coming from the PAC-Bayesian machine learning theory, arXiv:1209.0824.

\end{thebibliography}
\end{document}